\documentclass[12pt]{article}
\usepackage{graphpap}
\usepackage{amssymb}

\def \qed {\hfill $\boxempty$}

\def \complclass {\sf}
\def \np {{\complclass NP}}
\def \p {{\complclass P}}

\def \ttt {\mbox{\sf T}}
\def \fff {\mbox{\sf F}}

\newtheorem{Theorem}{Theorem}
\newtheorem{lem}[Theorem]{Lemma}
\newtheorem{defi}{Definition}
\newtheorem{crl}[Theorem]{Corollary}
\newtheorem{prp}[Theorem]{Proposition}
\newtheorem{rmk}[Theorem]{Remark}
\newtheorem{xmp}{Example}
\newtheorem{clm}{Claim}
\newtheorem{op}[Theorem]{Problem}

\def \komm {\LARGE\sf}

\def \bp {\begin{prp} \ }
\def \ep {\end{prp}}
\def \bc {\begin{crl} \ }
\def \ec {\end{crl}}
\def \thm {\begin{Theorem} \ }
\def \ethm {\end{Theorem}}
\def \bl {\begin{lem} \ }
\def \el {\end{lem}}
\def \bd {\begin{defi} \ \rm }
\def \ed {\end{defi}}
\def \brm {\begin{rmk} \ }
\def \erm {\end{rmk}}
\def \bxm {\begin{xmp} \ \rm }
\def \exm {\end{xmp}}
\def \bcm {\begin{clm} \ }
\def \ecm {\end{clm}}
\def \bop {\begin{op} \ }
\def \eop {\end{op}}
\def \nmr {\begin{enumerate}}
\def \enmr {\end{enumerate}}
\def \tmz {\begin{itemize}}
\def \etmz {\end{itemize}}

\def \nin {\noindent}
\def \bsk {\bigskip}
\def \ssk {\smallskip}
\def \msk {\medskip}
\def \pf {\nin{\bf Proof } \ }
\def \qed {\hfill $\Box$}
\def \enn {\mathbb{N}}
\def \eps {\varepsilon}
\def \vp {\varphi}
\def \aaa {\alpha}

\def \es {\emptyset}
\def \smin {\setminus}
\def \ssq {\subseteq}
\def \sst {\subset}

\def \NP {{\sf NP}}

\def \UU {\overline{\chi}}
\def \dec {\mbox{\rm dec}}

\def\cO{{\cal O}}
\def\cH{{\cal H}}
\def\cS{{\cal S}}

\def\cA{{\cal A}}
\def\cC{{\cal C}}
\def\cD{{\cal D}}

\def\cE{{\cal E}}

\def\cI{{\cal I}}

\begin{document}

\title{\vskip-1.5cm { ~ }
Approximability of the upper chromatic number
  of hypergraphs\thanks{
  ~Research supported in part by the Hungarian Scientific Research
  Fund, OTKA grant T-81493, and
 by the European Union and Hungary, co-financed
 by the European Social Fund through the project T\'AMOP-4.2.2.C-11/1/KONV-2012-0004 -- National Research Center
 for Development and Market Introduction of Advanced Information and Communication
 Technologies.}}
\author{Csilla Bujt\'as~$^{1}$ \ \qquad
   \vspace{2ex}
        Zsolt Tuza~$^{1,2}$
\\
\normalsize $^1$~Department of Computer Science
 and \vspace{-1mm} Systems Technology \\
 \normalsize University of  Pannonia,
 Veszpr\'em,
 Hungary \\
\normalsize $^2$~Alfr\'ed R\'enyi Institute of Mathematics\\
\normalsize        Hungarian Academy of Sciences,
 Budapest,
 \vspace{2mm} Hungary
 }
\date{\empty}
\maketitle

\begin{abstract}

A C-coloring of a hypergraph $\cH=(X,\cE)$ is a vertex coloring
$\vp:X\to\enn$ such that each edge $E\in\cE$ has at least two
vertices with a common color. The related parameter $\UU (\cH)$,
called the upper chromatic number of $\cH$, is the maximum number of
colors   can be used in a C-coloring of $\cH$. A hypertree is a
hypergraph which has a host tree $T$ such that each edge $E \in \cE$
induces a connected subgraph in $T$.
 Notations $n$ and $m$ stand for the number of vertices and edges, respectively, in a
generic input hypergraph.

 We establish guaranteed  polynomial-time approximation
  ratios   for the difference $n-\overline{\chi}({\cal H})$,
  which is $2+2 \ln (2m)$ on hypergraphs in general, and $1+ \ln m$
  on hypertrees.
  The latter ratio is essentially tight as we show that
  $n-\overline{\chi}({\cal H})$ cannot be approximated within
  $(1-\epsilon) \ln m$ on hypertrees
  (unless {\sf  NP}\,$\subseteq$\,{\sf DTIME}$(n^{\cO(log\;log\;
  n)})$).
   Furthermore, $\overline{\chi}({\cal H})$ does not have
  ${\cal O}(n^{1-\epsilon})$-approximation
    and
   cannot be approximated within additive error
  $o(n)$ on the class of hypertrees
  (unless  ${\sf P}={\sf NP}$).

\bigskip


\bigskip

\nin \textbf{AMS 2000 Subject Classification:}
 05C15, 
 05C65,  
 05B40, 
 68Q17 
\end{abstract}

\section{Introduction}

In this paper we study a hypergraph coloring invariant, termed upper
chromatic number and denoted by $\UU(\cH)$,   which was first
introduced by Berge (cf.~\cite{B})
 in the early 1970's
   and later independently by several further authors \cite{ABN, Vol2}
from different motivations.
 The
present work is the very first one concerning approximation
algorithms on it.

We also consider the complementary problem of approximating
 the difference $n-\UU$, the number of vertices minus the
 upper chromatic number.
One of our main tools to prove a guaranteed upper bound
  on it is an approximation ratio established for the
  2-transversal number of hypergraphs. As problems of this type are
  of  interest in their own right, we also prove an approximation ratio
  in general for the minimum size of  multiple transversals,
  i.e., sets of vertices intersecting each edge in a prescribed number
  of vertices at least.
%
 Earlier results allowed to select a vertex into the set several times;
  we prove bounds for the more restricted scenario where the set does
  not include any vertex more than once.

\subsection{Notation and terminology}

A \emph{hypergraph} $\cH=(X, \cE)$ is a set system, where $X$
denotes the set of vertices and each edge $E_i\in \cE$ is a nonempty
subset of $X$. Here we also assume that for each edge $E_i$ the
inequality $|E_i|\ge 2$ holds, moreover we use the standard
notations $|X|=n$ and $|\cE|=m$. A hypergraph $\cH$ is said to be
\emph{$r$-uniform} if $|E_i|=r$ for each $E_i \in \cE$.

 We shall also consider hypergraphs with restricted structure, where
 some kind  of host graphs are assumed. A hypergraph $\cH=(X,\cE)$
 admits a \emph{host graph} $G=(X,E)$ if each edge $E_i \in \cE$ induces a
 connected subgraph in $G$.
 The edges of the host graph $G$ will be  referred to as
 \emph{lines}.
 Particularly, $\cH$ is called
 \emph{hypertree}
 or  \emph{hyperstar} if it admits a host graph which is a tree or a star,
 respectively. Note that under our condition, which forbids edges of
 size 1,
  $\cH$ is a hyperstar if and only if there exists a
 fixed
 vertex $c^*\in X$ (termed the center of the hyperstar) contained in
 each edge of $\cH$.

 \bsk

 A \emph{C-coloring} of $\cH$ is an assignment $\vp:X\to\enn$ such that
 each edge $E\in\cE$ has at least two
vertices of a common color (that is, with the same image).
   The \emph{upper chromatic number} $\UU(\cH)$ of $\cH$ is the maximum
    number of colors that can be used in a C-coloring of $\cH$.
  We note that in the literature the value
 $\UU(\cH)+1$ is also called the `cochromatic
 number' or `heterochromatic number' of $\cH$
 with the  terminology of Berge \cite[p.~151]{B}
 and Arocha \emph{et al.}~\cite{ABN},
 respectively.
 A C-coloring $\vp$ with  $|\vp(X)|=\UU(\cH)$ colors
  will be referred to as an \emph{optimal coloring} of $\cH$.
  The \emph{decrement} of $\cH=(X,\cE)$,   introduced in
  \cite{proj-plane},
   is defined as $\dec(\cH)=n-\UU(\cH)$. Similarly, the decrement
   of a C-coloring $\vp:X\to\enn$ is meant as
   $\dec(\vp)=|X|-|\vp(X)|$.
   For results on C-coloring see the recent survey \cite{BT-JGeom}.

 \bsk

  A   \emph{transversal}  (also called hitting set  or   vertex cover)
  is a subset $T \subseteq X$ which meets each edge of $\cH=(X, \cE)$,
 and the minimum cardinality  of a transversal  is the
 \emph{transversal number} $\tau(\cH)$ of the hypergraph.
 An \emph{independent set} (or stable set) is a vertex set $I \subseteq
 X$, which contains no edge of $\cH$ entirely. The maximum
 size of an independent set in $\cH$ is the \emph{independence number} (or
 stability number) $\aaa(\cH)$. It is immediate from the definitions
 that the complement of a transversal is an independent set and vice versa,
 so the Gallai-type equality $\tau(\cH)+\aaa(\cH)=n$ holds for
 each hypergraph.
 Remark that selecting one vertex from each color class of a C-coloring
   yields an independent set, therefore $\UU(\cH)\le\aaa(\cH)$ and,
   equivalently, $\dec(\cH) \ge \tau (\cH)$.

 More generally, a \emph{$k$-transversal} is a set $T\subseteq X$ such that
 $|E_i \cap T|\ge k$ for every $E_i \in \cE$. A 2-transversal is
 sometimes called double transversal or strong transversal, and its
 minimum size is the \emph{2-transversal number} $\tau_2(\cH)$ of the
 hypergraph.

 \bsk

  For an optimization problem and a constant $c>1$, an algorithm $\cA$
is
 called a \emph{$c$-approximation algorithm} if, for every feasible
 instance $\cI$ of the problem,
%
\tmz
 \item
  if the value has to be minimized, then $\cA$ delivers a
  solution of value at most $c\cdot Opt(\cI)$;
 \item
  if the value has to be maximized, then $\cA$ delivers a
  solution of value at least $Opt(\cI)/c$.
\etmz

 Throughout this paper, an approximation algorithm is always meant
 to be one with polynomial running time on every instance of the problem.
  We say that a value has  guaranteed approximation ratio $c$
  if it has  a $c$-approximation algorithm. In the other case, when
  no $c$-approximation algorithm exists, we say that the value
  cannot be approximated within   ratio $c$.
  For a function $f(n,m)$, an $f(n,m)$-approximation algorithm
   and the related notions can be defined
 similarly.
 A polynomial-time approximation scheme, abbreviated as PTAS, means an
  algorithm for every fixed $\eps > 0$ which
   is a $(1+\eps)$-approximation and
  whose running time is a polynomial function of the input size
  (but any function of $1/\eps$ may occur in the exponent).

  \bsk

For further terminology and facts we refer to \cite{B,BM,Vaz}
 in the theory of graphs, hypergraphs,
 and algorithms, respectively.
 The notations $\ln x$ and $\log x$ stand for the
 natural logarithm and for the logarithm in base 2, respectively.

\subsection{Approximability results on multiple transversals}

 The  transversal number $\tau(\cH)$ of a hypergraph can be approximated within
 ratio $(1+\ln m)$  by the classical greedy algorithm (see e.g. \cite{Vaz}). On the other
 hand,   Feige \cite{Fei} proved that $\tau(\cH)$ cannot be
 approximated within $(1-\epsilon) \ln m$ for any constant
 $0<\epsilon <1$,  unless {\sf  NP}\,$\subseteq$\,{\sf DTIME}$(n^{\cO(\log \log  n)})$.
 As relates to the $k$-transversal number, in \cite{Vaz} a $(1+ \ln m)$-approximation
 is stated under the less restricted setting which allows multiple selection of
 vertices in the $k$-transversal.
 In the context of coloring, however, we cannot allow repetitions of vertices.
  For this more restricted case, when
  the $k$-transversal consists of pairwise different vertices, we
 prove a guaranteed approximation ratio $(1+ \ln (km))$.

 In fact we consider a more general problem, where the required minimum size
 of the intersection $E_i \cap T$ can be prescribed independently
 for each $E_i \in \cE$.
 \thm
 \label{multiple}
  Given a hypergraph\/ $\cH=(X,\cE)$ with\/ $m$ edges\/
   $E_1,\dots,E_m$ and positive integers\/ $w_1,\dots,w_m$
   associated with the edges, the minimum cardinality of a set\/
   $S\sst X$ satisfying\/ $|S\cap E_i|\ge w_i$ for all\/
   $1\le i\le m$ can be approximated within\/
   $\sum_{i=1}^{W} 1/i < 1+\ln W$,
    where\/ $W= \sum_{i=1}^m w_i$.
\ethm

 This result, proved in the next section,
  implies a guaranteed approximation ratio $(1+ \ln 2m)$ for
 $\tau_2(\cH)$.

 \subsection{Approximability results on  the upper chromatic number}

 The problem of determining the upper chromatic number
 is \NP-hard, already on the class of 3-uniform hyperstars.
 On the other hand,  the problems of determining $\overline{\chi}({\cal H})$ and finding a
  $\overline{\chi}({\cal H})$-coloring are fixed-parameter tractable in terms of
  maximum vertex degree on the class of hypertrees \cite{BT-cejor}.

  \bsk

  A notion closely related to our present subject was introduced by
  Voloshin \cite{Vol93, Vol2} in 1993. A \emph{mixed hypergraph} is a
  triple $\cH =(X, \cC, \cD)$  with
 two families of subsets called $\cC$-edges and $\cD$-edges.
     By definition, a coloring of a mixed hypergraph is an
     assignment $\vp:X\to\enn$ such that each $\cC$-edge has two
     vertices of a common color and each $\cD$-edge has two vertices
     of distinct colors. Then, the minimum and the maximum possible
     number of colors, that can occur in a coloring of $\cH$, is termed
     the lower and the upper chromatic number of $\cH$ and denoted
     by $\chi(\cH)$ and $\UU(\cH)$, respectively. For detailed
     results on mixed hypergraphs we refer to the monograph \cite{Volmon}.
     Clearly, the \mbox{C-colorings} of a hypergraph $\cH=(X, \cE)$ are in
     one-to-one correspondence with the colorings of the mixed hypergraph
      $\cH'=(X, \cE, \es)$, and also $\UU(\cH)= \UU(\cH')$ holds.

 \bsk
 The following results are known on the approximation of the upper
 chromatic number of mixed hypergraphs:
 \tmz
  \item
  For mixed hypergraphs of maximum degree 2, the upper chromatic
   number has a linear-time $\frac{5}{3}$-approximation and an
   $O(m^3+n)$-time \mbox{$\frac{3}{2}$-approximation.}
    \cite[Theorem 14 and Theorem 15]{KKV-degree}
 \item
   There is no PTAS for the upper
   chromatic number of mixed hypergraphs of maximum degree 2, unless \p\,$=$\,\np.
    \cite[Theorem 20]{KKV-degree}
  \item
  There is no  $o(n)$-approximation algorithm
   for the upper chromatic number of mixed hypergraphs,
   unless \p\,$=$\,\np.
    \cite[Corollary 5]{K-spect}
 \etmz
All these results assume the presence of $\cD$-edges in the input
mixed hypergraph. In this paper we investigate how hard it is to
estimate $\UU$ for C-colorings of  hypergraphs.

\bsk

 On the positive side, we prove a guaranteed approximation ratio for
 the decrement  of hypergraphs in general, furthermore we establish a better
 ratio on the class of hypertrees.

\thm \label{appr-gen}
 The value of\/ $\dec(\cH)$ is\/
 $(2+2\ln (2m))$-approximable on the class of all hypergraphs.
\ethm

\thm \label{appr-htree}
 The value of\/ $\dec(\cH)$ is\/
 $(1+\ln m)$-approximable on the class of all hypertrees.
\ethm

These theorems are essentially best possible concerning the ratio of
approximation, moreover the upper chromatic number turns out to be
inherently non-approximable already on hypertrees with rather
restricted host trees, as shown by the next result.

\thm
    \label{ratio}
   \tmz
 \item[$(i)$]
  For every\/ $\epsilon > 0$, \enskip $\dec(\cH)$ cannot be approximated
 within\/ $(1-\epsilon)\ln m$ on the class of hyperstars,
  unless\/ {\sf  NP}\,$\subseteq$\,{\sf DTIME}$(n^{\cO(\log \log  n)})$.
    \item[$(ii)$]
   For every\/ $\epsilon > 0$, \enskip $\UU(\cH)$  cannot be approximated
  within\/   $n^{1-\epsilon}$  on the class of\/ $3$-uniform hyperstars,
  unless\/ {\sf P}\,$=$\,{\sf NP}.
  \etmz
  \ethm

As regards the \emph{difference} between a solution determined by
a polynomial-time algorithm and the optimum value, the situation
is even worse.

\thm
    \label{additive}
 Unless\/ \p\,$=$\,\np, neither of the following values
  can be approximated within additive error\/ $o(n)$ for
  hypertrees of edge size at most 7\,:\,

\msk

  $\UU(\cH)$, \ $\dec(\cH)$, \
  $\aaa(\cH)-\UU(\cH)$, \ $\tau(\cH)-\dec(\cH)$, \
  $\dec(\cH)-\tau_2(\cH)/2$.
 \ethm

The relevance of the last quantity occurs in the context of
Proposition \ref{decr-tau2} of Section \ref{decr-transv}.

We prove the positive results with guaranteed approximation ratio
 in Section 3, and the negative non-approximability results
 in Section 4.

\subsection{Lemmas on connected colorings of hypertrees}

Suppose that $\cH$ is a hypergraph over a host graph $G$, and $\vp$ is a
 C-coloring of $\cH$.
We say that $\vp$ is a \emph{connected coloring} if each color class
  of $\vp$ induces a connected subgraph of $G$.
 We will use the following two lemmas concerning connected C-colorings
 of hypertrees,
   both   established in \cite{BT-cejor}.
   A line $uv$ of the host tree $G$
   is termed \emph{monochromatic line} for a C-coloring $\vp$ if $\vp(u)=\vp(v)$.
    \bl(\cite[Proposition 2]{BT-cejor})
     \label{conn}
      If a hypertree admits a C-coloring with\/ $k$ colors, then it
      also has a connected C-coloring with\/ $k$ colors over any fixed
      host tree.
    \el

    \bl(\cite[Proposition 3]{BT-cejor})
    \label{mono-lines}
    If\/ $\vp$ is a connected C-coloring of a hypertree\/ $\cH$ over a
    fixed host tree\/ $G$, then the decrement of\/ $\vp$ equals the
    number of monochromatic lines in\/ $G$.
    \el

\section{Multiple transversals}

In this section, we describe a variation of the classical greedy
algorithm, with the goal to produce a multiple transversal  with pairwise
different elements.  Analyzing the greedy selection we will prove
Theorem \ref{multiple}. We recall its statement.
  \bsk

 \noindent \textbf{Theorem~\ref{multiple}}.
 \emph{Given a hypergraph\/ $\cH=(X,\cE)$ with\/ $m$ edges\/
   $E_1,\dots,E_m$ and positive integers\/ $w_1,\dots,w_m$
   associated with its edges, the minimum cardinality of a set\/
   $S\sst X$ satisfying\/ $|S\cap E_i|\ge w_i$ for all\/
   $1\le i\le m$ can be approximated within\/
   $\sum_{i=1}^{W} 1/i < 1+\ln W$,
    where\/ $W= \sum_{i=1}^m w_i$.}

    \msk
    \pf
 Denote by $\cS$ the collection of all feasible solutions,
  that are the sets $S\sst X$ such that
  $|S\cap E_i|\ge w_i$ holds for all $i=1,\dots,m$.
 By definition, the optimum of the problem is the integer
    $$M:=\min_{S\in\cS} |S| .$$
 We will show that the greedy  selection always yields an
$S^*\in\cS$ with
 $$|S^*| \le M \cdot \left(1 + 1/2 + \dots + 1/W\right).$$
 To prove this, for any $Y\sst X$ and any $1\le i\le m$ we define
 $$w_{i,Y} := \max \left(0, \, w_i - |E_i\cap Y|\right)$$
 which means the reduced number of elements to be picked
 further from $E_i$, once the set $Y$ has already been selected.
Moreover, to any vertex $x\in X\smin Y$ we associate its usefulness
 $$u_{x,Y} := | \{E_i \mid x\in E_i, \ w_{i,Y} > 0 \}|.$$
The greedy algorithm then starts with $Y_0=\es$ and updates
 $Y_k := Y_{k-1}\cup\{x_k\}$ where $x_k\in X\smin Y_{k-1}$ has maximum
  usefulness among all values $u_{x,Y_{k-1}}$ in the set $X\smin Y_{k-1}$,
   as long as this maximum is positive.
Reaching $u_{x,Y_t}=0$ for all $x\in X\smin Y_t$ (for some $t$),
 we set $S^* := Y_t$; we will prove that this $S^*$ satisfies the requirements.

It is clear by the definition of $u_{x,Y}$ that $S^*$ meets
 each $E_i$ in at least $w_i$ elements, i.e.\ $S^*\in\cS$.
We need to prove that $S^*$ is sufficiently small.
 For this, consider the following auxiliary set of cardinality $W$:
 $$Z := \{ z(i,j) \mid 1\le i\le m, \ 1\le j\le w_i \}.$$
At the moment when $Y_k$ is constructed by adjoining
 an element $x_k$ to $Y_{k-1}$, we assign weight $1/u_{x,Y_{k-1}}$ to
  all elements $z(i,w_{i,Y_{k-1}})$ such that
  $x_k\in E_i$ and $w_{i,Y_{k-1}}>0$.
Note that $w_{i,Y_{k}}=w_{i,Y_{k-1}}-1$ will hold after the
 selection of $x_k$.
Moreover, total weight 1 is assigned in each step, hence the overall
 weight after finishing the algorithm is exactly $|S^*|$.
We put the elements $z(i,j)$ in a sequence $Z^*=(z_1,z_2,\dots,z_W)$
 such that the elements of $Z$ occur in the order as they are weighted
 (i.e., those for $x_1$ first in any order,
  then the elements weighted for $x_2$, and so on).

Just before the selection of $x_k$, the number of elements $z(i,j)$
 to which a weight has been assigned is precisely
 $m_{k-1} := \sum_{\ell=1}^{k-1} u_{x_\ell,Y_{\ell-1}}.$
  We are going to prove that $u_{x_k,Y_{k-1}} \ge (W-m_{k-1})/M$.
Assuming that this has already been shown, it follows that
 each $z_q$ in $Z^*$ has weight at most $M/(W+1-q)$ and consequently
 $|S^*| \le M \cdot \left(1 + 1/2 + \dots + 1/W\right)$ as required.

 Let now $S_0\in\cS$ be any fixed optimal solution.
Consider the bipartite incidence graph $B$ between the sets $E_i$
and the
 elements of $S_0$.
That is, the first vertex class of $B$ has $m$ elements
$a_1,\dots,a_m$
 representing the sets $E_1,\dots,E_m$ while the second vertex class
 consists of the elements of $S_0$; we denote the latter vertices by
 $b_1,\dots,b_M$.
There is an edge joining $a_i$ with $b_j$ if and only if $b_j\in
E_i$.

 Since $S_0\in \cS$, each $a_i$ has degree at least $w_i$.
Moreover, considering the moment just before $x_k$ is selected, if
we
 remove the vertices of $S_0\cap Y_{k-1}$, in the remaining subgraph
 still each $a_i$ has degree at least $w_{i,Y_{k-1}}$.
We take a subgraph $B'$ of this $B-Y_{k-1}$ (possibly $B$ itself if
 $Y_{k-1}\cap S_0=\es$) such that each $a_i$ has degree
  \emph{exactly} $w_{i,Y_{k-1}}$.
The number of edges in $B'$ is then equal to $W-m_{k-1}$;
 hence, some $b_j$ has degree at least $(W-m_{k-1})/M$.
It follows that this $b_j$ has usefulness at least $(W-m_{k-1})/M$
 at the moment when $x_k$ is selected; but $x_k$ is chosen to have
 maximum usefulness, hence $u_{x_k,Y_{k-1}} \ge (W-m_{k-1})/M$.
This completes the proof.
  \qed

\bc \label{k-transv}
 For each positive integer\/ $k$, the\/ $k$-transversal number\/ $\tau_k$
 has a\/ $(1+\ln (km))$-approximation on the class of all
 hypergraphs.
\ec

\section{Guaranteed approximation ratios for the decrement}
In this section we establish a connection between the parameters
$\dec(\cH)$   and $\tau_2(\cH)$, and then we prove our positive
results stated in Theorems~\ref{appr-gen} and
\ref{appr-htree}.

\subsection{Decrement vs.\ 2-transversal number}
   \label{decr-transv}

First,  we give an inequality valid for all hypergraphs without
 any structural restrictions and then, using this relation, we prove
 Theorem~\ref{appr-gen}.

\bp   \label{decr-tau2}
 For every hypergraph\/ $\cH$ we have\/
  $\tau_2(\cH)/2\le\dec(\cH)\le\tau_2(\cH)-1$, and
  both bounds are tight.
 In particular,\/
  $\tau_2(\cH)$ is a 2-approximation for\/ $\dec(\cH)$.
\ep

 \pf
{\sl Lower bound:}\quad If $\UU(\cH) \le n/2$, then $\dec(\cH)\ge
n/2 \ge \tau_2(\cH)/2$ automatically holds. If $\UU(\cH) > n/2$,
then every $\UU$-coloring contains at least $2\UU(\cH)-n$ singleton
color classes, therefore the total size of non-singleton classes is
at most $n-(2\UU(\cH)-n)= 2(n-\UU(\cH))$. Since the union of the
latter meets all edges at least twice, we obtain
$2\dec(\cH)\ge\tau_2(\cH)$.

\ssk

{\sl Upper bound:}\quad If $S$ is a 2-transversal set of cardinality
$\tau_2(\cH)$, we can assign the same color to the entire $S$ and a
new dedicated color to each $x\in X\smin S$. This is a C-coloring
with $n-|S|+1$ colors and with decrement $\tau_2(\cH)-1$.

\ssk

{\sl Tightness:}\quad The simplest example for equality in the upper
bound is the hypergraph in which the vertex set is the only edge,
i.e.\ $\cH=(X,\{X\})$. Many more examples can be given. For
instance, we can specify a proper subset $S\sst X$ with $|S|\ge 2$,
and take all triples
$E\sst X$ such that $|E\cap S|=2$ and $|E\smin S|=1$. If $|S|\le
n-2$, then $S$ is the unique smallest 2-transversal set, and every
C-coloring with more than two colors makes $S$ monochromatic,
hence the unique $\UU$-coloring uses $n-|S|+1$ colors.

For the lower bound, we assume that $n=3k+1$. Let $X=\{1,2, \dots,
3k+1\}$ and
\begin{eqnarray}
  \cE&=&\{\{3r+1, 3r+2, 3r+3\}\mid 0\le r\le k-1\} \nonumber \\
   & &\cup~ \{\{3r+2, 3r+3, 3r+4\}\mid 0\le r\le k-1\}\} \nonumber
\end{eqnarray}
Then $\tau_2(\cH)=2k$ because the $k$ edges in the first line are
mutually disjoint and hence need at least $2k$ vertices in any
2-transversal set, while the $2k$-element set $\{3r+2\mid 0\le r\le
k-1\}\cup \{3r+3\mid 0\le r\le k-1\}$ meets all edges twice. On the
other hand, there exists a unique C-coloring with decrement $k$,
obtained by making $\{3r+2,3r+3\}$ a monochromatic pair for
$r=0,1,\dots,k-1$ and putting any other vertex in a singleton color
class. This verifies equality in the lower bound.
  \qed

\bsk

 \nin Now, we are ready to prove Theorem~\ref{appr-gen}. Let us recall its
 statement.

 \bsk
 \noindent \textbf{Theorem~\ref{appr-gen}}.
 \emph{The value of\/ $\dec(\cH)$ is\/
 $(2+2\ln (2m))$-approximable on the class of all hypergraphs.}
 \msk

 \pf By Corollary \ref{k-transv}, we have a  $(1+ \ln (2m))$-approximation
 algorithm $\cA$ for $\tau_2$.
  Hence, given a  hypergraph $\cH=(X,\cE)$, the algorithm $\cA$
 outputs a 2-transversal $T$ of size at most $(1+ \ln (2m))\tau_2(\cH)$.
  Then, assign color 1 to every $x\in T$, and color the $n-|T|$ vertices
  in $X\setminus T$ pairwise differently with colors $2,3,\dots,
  n-|T|+1$. As each edge $E_i\in \cE$ contains at least two vertices
  of color 1, this results in a C-coloring $\vp$ with decrement
  satisfying
  $$\dec(\vp) = |T|-1 \le (1+ \ln (2m))\tau_2(\cH) -1 < 2(1+ \ln
  (2m))\dec(\cH),$$
  where the last inequality follows from
  Proposition~\ref{decr-tau2}.
  Therefore, algorithm $\cA$ together with the simple construction
  of coloring $\vp$ is a $(2+2\ln 2m)$-approximation for $\dec(\cH)$.
  \qed

  \subsection{Guaranteed approximation ratio on hypertrees}

   In this short subsection we prove  Theorem~\ref{appr-htree}. We recall its
   statement.

    \bsk
 \noindent \textbf{Theorem~\ref{appr-htree}}.
 \emph{The value of\/ $\dec(\cH)$ is\/
 $(1+\ln m)$-approximable on the class of all hypertrees.}

   \msk
   \pf Given a hypertree $\cH=(X, \cE)$ and $G=(X,L)$ which is a
   host tree of $\cH$, construct the auxiliary hypergraph $\cH^*=(L^*,
   \cE^*)$ such that each vertex $l_i^* \in L^*$ represents a line
   $l_i$ of the host tree, moreover each edge $E_i^* \in \cE^*$ of
   the auxiliary hypergraph corresponds to the edge $E_i  \in \cE$
   in the following way:
   $$E_i^*=\{l_j^* \mid l_j \subseteq E_i\}.$$

   Now, consider any connected C-coloring  $\vp$ of $\cH$.
    This coloring determines the set $S
   \subseteq L$ of monochromatic lines in the host tree, moreover the
   corresponding vertex  set $S^* \subseteq L^*$ in $\cH^*$.
   By Lemma~\ref{mono-lines}, $\dec(\vp)=|S|=|S^*|$.
      As $\vp$ is a connected C-coloring, each edge of $\cH$ contains
   a monochromatic line and, consequently, $S^*$ is a transversal of size
   $\dec(\vp)$ in    $\cH^*$.
    Similarly, in the opposite direction, if a transversal
   $T^*$ of $\cH^*$ is given and  the corresponding line-set is $T$
   in the host tree, then every edge $E_i$ of $\cH$ contains two vertices,
   say $u$ and $v$, such that the line $uv$ is contained in $T$.
   Then, the vertex coloring $\phi$, whose color classes correspond
   to the components of $(X, T)$, is a connected C-coloring of
   $\cH$, and in addition $\dec(\phi)=|T|=|T^*|$ holds.

   By Lemma~\ref{conn}, $\cH$ has a connected C-coloring  $\vp$ with
   $\dec(\vp)=\dec(\cH)$, therefore the  correspondence above implies
   $\dec(\cH)= \tau(\cH^*)$.

   As $\cH^*$ can be constructed in polynomial time from the hypertree $\cH$,
   and since a transversal $T^*$ of size at most $(1+ \ln m)\tau(\cH^*)$ can be
   obtained by  greedy selection, a C-coloring $\phi$ of $\cH$ with
   $$\dec(\phi)=|T^*| \le (1+ \ln m)\tau(\cH^*)=(1+ \ln m)\dec(\cH)$$
   can also be constructed in polynomial time.
   This yields  a guaranteed approximation ratio $(1+ \ln m)$ for
   the decrement on the class of hypertrees.
   \qed

\section{Approximation hardness}

 The bulk of this section is devoted to the proof of Theorem~ \ref{additive}
 on non-approximability for hypertrees.
  Then, we prove a lemma concerning parameters $\UU(\cH)$ and $\dec(\cH)$ of
  hyperstars. The section is closed with the proof of
  Theorem~\ref{ratio} and with some remarks.

\subsection{Additive linear error}

 Our goal in this subsection is to prove Theorem \ref{additive}.
  This needs the following construction,
  which was introduced in \cite{perf-htree}.
  (We note that a similar  construction was given already in \cite{KKPV}.)

\paragraph{Construction of $\cH(\Phi)$.}

Let $\Phi= C_1 \wedge \cdots \wedge C_m$ be an instance of 3-SAT,
with $m$ clauses of size 3 over the set $\{x_1,\dots,x_n\}$ of $n$
variables, such that the three literals in each clause $C_j$ of
$\Phi$ correspond to exactly three distinct variables. We
construct the hypertree $\cH=\cH(\Phi)$ with the set
  $$
  X = \{ c^* \}
   \cup \{ x'_i,\, t_i,\, f_i \mid 1 \le i \le n \}
 $$
of $3n+1$ vertices, where the  vertices $x'_i,t_i,f_i$ correspond to
variable $x_i$. First, we define the host tree $T=(X,E)$ with vertex
set $X$ and line-set
$$
  E = \{ c^*x_i',\, x_i't_i,\, x_i'f_i \mid 1\leq i\leq n \}.
$$

Hypergraph $\cH$ will have 3-element ``variable-edges''
$H_i=\{x_i',t_i,f_i\}$ for $i=1,\dots,n$, and 7-element
``clause-edges'' $F_j$ representing clause $C_j$ for
$j=1,\dots,m$. All the latter contain $c^*$ and six further
vertices, two for each literal of $C_j$\,:
  \tmz
   \item If $C_j$ contains the positive literal $x_i$, then  $F_j$
  contains $x_i'$ and $t_i$.
  \item If $C_j$ contains the negative literal $\neg x_i$, then  $F_j$
  contains $x_i'$ and $f_i$.
  \etmz
Since $H_1,\dots,H_n$ are disjoint edges, it is clear that
$\dec(\cH)\ge n$ and $\UU(\cH)\le 2n+1$. We shall see later that
equality holds if and only if $\Phi$ is satisfiable. In addition,
since $x_1',\dots,x_n'$ is a transversal set of $\cH$, the
equalities $\tau(\cH)=n$ and $\aaa(\cH)=2n+1$ are valid for all
$\Phi$, no matter whether satisfiable or not. Also, $\tau_2(\cH)=2n$
for all $\Phi$.

\paragraph{Optimal colorings of $\cH$.}

By Lemma \ref{conn}, we may restrict our attention to colorings
where each color class is a subtree in $T$. This makes a coloring
irrelevant if it 2-colors a variable-edge in such a way that
$\{t_i,f_i\}$ is monochromatic but $x_i'$ has a different color.
Hence, at least one of the lines $x_i't_i$ and $x_i'f_i$ is
monochromatic (maybe both) for each $i$. Moreover, we may assume the
following further simplification: there is no monochromatic line
 $c^*x_i'$. Indeed, if the entire $H_i$ is monochromatic, then we
 would lose a color by making the line $c^*x_i'$ monochromatic. On
 the other hand, if say the monochromatic pair inside $H_i$ is
$x_i't_i$, then every clause-edge $F_j$ containing $c^*x_i'$ but
avoiding $t_i$ also contains the line $x_i'f_i$, therefore we get a
coloring with the same number of colors if we assume that $x_i'f_i$
is monochromatic instead of $c^*x_i'$. Summarizing, we search an
optimal coloring $\vp:X\to\enn$ with the following properties for
all $i=1,\dots,n$\,:
 \tmz
  \item $\vp(c^*)\ne\vp(x_i')$
  \item $\vp(x_i')=\vp(t_i)$ or $\vp(x_i')=\vp(f_i)$
 \etmz
In the rest of the proof we assume that all vertex colorings
occurring satisfy these conditions.

\paragraph{Truth assignments.}

Given a coloring $\vp$, we interpret it in the following way for
truth assignment and clause deletion:
 \tmz
  \item If $H_i$ is monochromatic, delete all clauses from $\Phi$
   which contain literal $x_i$ or $\neg x_i$.
  \item Otherwise, assign truth value
    $x_i\mapsto\ttt$ if $\vp(x_i')=\vp(t_i)$,
    and $x_i\mapsto\fff$ if $\vp(x_i')=\vp(f_i)$.
 \etmz
It follows from the definition of $\cH(\Phi)$ that this truth
assignment satisfies the modified formula after deletion if and
only if $\vp$ properly colors all edges of $\cH$.

Also conversely, if $\Phi'$ is obtained from $\Phi$ by deleting all
clauses which contain $x_i$ or $\neg x_i$ for a specified index set
$I\ssq\{1,\dots,n\}$, then a truth assignment $a:\{x_i \mid i\in
\{1,\dots,n\}\smin I \}\to\{\ttt,\fff\}$ satisfies $\Phi'$ if and
only if the following specifications for the monochromatic lines
yield a proper coloring $\vp$ of $\cH$\,:
 \tmz
  \item If $i\in I$, then $\vp(x_i')=\vp(t_i)=\vp(f_i)$.
  \item Otherwise, let $\vp(x_i')=\vp(t_i)$ if $a(x_i)=\ttt$,
   and $\vp(x_i')=\vp(f_i)$ if $a(x_i)=\fff$.
 \etmz
The observations above imply the following statement:

\bl
 For any instance\/ $\Phi$ of\/ {\rm 3-SAT}, the value of\/
  $\dec(\cH(\Phi))$ is equal
  to the minimum number of variables whose deletion from\/
  $\Phi$ makes the formula satisfiable.
   \qed
\el

To complete our preparations for the proof of the theorem, let us
quote an earlier result on formulas in which every positive and
negative literal occurs in at most four clauses. The problem {\sc
Max 3Sat$(4,\overline 4)$} requires to maximize the number of
satisfied clauses in such formulas. The following assertion states
that this optimization problem is hard to approximate, even when
the input is restricted to satisfiable formulas.

\bl   (\cite[Corollary 5]{bounded-sat}) \
   \label{sat-bounded}
 Satisfiable\/ {\sc Max 3Sat$(4,\overline 4)$} has no PTAS,
 unless\/ \p\,$=$\,\np.
\el

  Now we are ready to verify Theorem~\ref{additive}, which states:
 \bsk

 \noindent \textbf{Theorem~\ref{additive}}.
  \emph{Unless\/ \p\,$=$\,\np, neither of the following values
  can be approximated within additive error\/ $o(n)$ for
  hypertrees of edge size at most 7\,:\,}

\msk

  $\UU(\cH)$, \ $\dec(\cH)$, \
  $\aaa(\cH)-\UU(\cH)$, \ $\tau(\cH)-\dec(\cH)$, \
  $\dec(\cH)-\tau_2(\cH)/2$.

  \msk
  \pf
  We apply reduction from Satisfiable {\sc Max 3Sat$(4,\overline 4)$}.
For each instance $\Phi$ of this problem, we construct the
hypergraph $\cH=\cH(\Phi)$. Since $\Phi$ is required to be
satisfied, no variables have to be deleted from it to admit a
satisfying truth assignment. This means precisely one monochromatic
line inside each variable-edge. Hence, the above observations
together with Lemma \ref{conn} imply that $\dec(\cH)=n$ and
$\UU(\cH)=2n+1$.

On the other hand, Lemma \ref{sat-bounded} implies the existence
of a constant $c>0$ such that it is \np-hard to find a truth
assignment that satisfies all but at most $cm$ clauses in a
satisfiable instance of {\sc Max 3Sat$(4,\overline 4)$} with $m$
clauses. Since each literal occurs in at most four clauses, this
may require the cancelation of at least $cm/8\ge c'n$ variables.
Thus, for the coloring $\vp$ determined by a polynomial-time
algorithm, $\dec(\vp)-\dec(\cH)=\Theta(n)$ may hold, and hence
also $\UU(\cH)-|\vp(X)|=\Theta(n)$.
 \qed

  \subsection{No efficient approximation on hyperstars}

 Proposition~\ref{decr-tau2}   established a relation between
 $\dec(\cH)$ and $\tau_2(\cH)$, valid for all hypergraphs.
  Here we show that for hyperstars there is a stronger
  correspondence between the parameters. After that,
     we prove Theorem \ref{ratio}
      which states
   non-approximability results on hyperstars.
   \bsk

 Given a \emph{hyperstar} $\cH=(X,\cE)$,  let
  us denote by $c^*$ the
center of the host star. Hence, $c^*\in E$ holds for all $E\in\cE$.
We shall use the following notations:
$$
  E^- = E \smin \{c^*\}, \quad
  \cE^- = \{E^-\mid E\in\cE\}, \quad
  \cH^- = (X\smin\{c^*\},\cE^-) .
$$

\bp   \label{star-decrem}
 If\/ $\cH$ is a hyperstar, then\/
  $\dec(\cH) = \tau(\cH^-) = \tau_2(\cH) -1$ and\/
  $\UU(\cH)=\aaa(\cH^-)+1$.
\ep

 \pf
If a 2-transversal set $S$ does not contain $c^*$, then we can
replace any $s\in S$ with $c^*$ and obtain another 2-transversal set
of the same cardinality. This implies $\tau(\cH^-) = \tau_2(\cH)
-1$.

Let us observe next that the equalities $\UU(\cH)=\aaa(\cH^-)+1$ and
$\dec(\cH) = \tau(\cH^-)$ are equivalent, due to the Gallai-type
equality for $\aaa+\tau$ in $\cH^-$.

Now, the particular case of Lemma \ref{conn} for hyperstars means
that there exists a $\UU$-coloring of $\cH$ such that all color
classes but that of $c^*$ are singletons. Those singletons form an
independent set in $\cH^-$, because the color of $c^*$ is repeated
inside each $E^-$. Thus, we necessarily have
$\UU(\cH)\le\aaa(\cH^-)+1$.

Conversely, if $S$ is a largest independent set in $\cH^-$, i.e.\
$|S| = \aaa(\cH^-) = |X|-1-\tau(\cH^-)$ and $E^-\smin S\ne\es$ for
all $E^-$, then making $X\smin S$ a color class creates a
monochromatic pair inside each $E\in\cE$ because the color of $c^*$
is repeated in each $E^-$. Hence, assigning a new private  color to
each $x\in S$ we obtain that $\UU(\cH)\ge\aaa(\cH^-)+1$,
consequently $\UU(\cH)=\aaa(\cH^-)+1$ and $\dec(\cH)=\tau(\cH^-)$.
 \qed

\bsk

 The following non-approximability results concerning $\UU(\cH)$ and $\dec(\cH)$
 are valid already on the class of
  hyperstars. We recall the statement of Theorem~\ref{ratio}.

\bsk

 \noindent \textbf{Theorem~\ref{ratio}}.
  \emph{\tmz
 \item[$(i)$]
 For every\/ $\epsilon > 0$, \enskip $\dec(\cH)$ cannot be approximated
 within\/ $(1-\epsilon)\ln m$ on the class of hyperstars,
  unless\/ {\sf  NP}\,$\subseteq$\,{\sf DTIME}$(n^{\cO(\log \log  n)})$.
  \item[$(ii)$]
  For every\/ $\epsilon > 0$, \enskip $\UU(\cH)$  cannot be approximated
  within\/   $n^{1-\epsilon}$  on the class of\/ $3$-uniform hyperstars, unless\/ {\sf
 P}\,$=$\,{\sf NP}.
  \etmz}
  \msk

\pf By Proposition \ref{star-decrem}, the equalities
$\UU(\cH)=\aaa(\cH^-)+1$ and $\dec(\cH)=\tau(\cH^-)$ hold whenever
$\cH$ is a hyperstar.

 \tmz
 \item[$(i)$]
  If $\cH$ is a generic hyperstar (with no
 restrictions on its edges), then $\cH^-$ is a generic hypergraph.
 Thus, approximating $\dec(\cH)$ on hyperstars is equivalent to
 pproximating $\tau(\cH^-)$ on hypergraphs, which is known to be
 intractable within ratio $(1-\eps)(\log m)$ unless {\sf
 NP}\,$\subseteq$ {\sf DTIME}$(n^{\cO(\log \log  n)})$, by the result
 of Feige \cite{Fei}.
    \item[$(ii)$]
If $\cH$ is a generic 3-uniform hyperstar, then $\cH^-$ is a generic
graph. Thus, approximating $\UU(\cH)$ on 3-uniform hyperstars is
equivalent to approximating $\aaa(\cH^-)+1$ on graphs, which is
known to be intractable within ratio $n^{1-\eps}$ unless \p\,=\,\np,
by the result of
Zuckerman \cite{Zuck}.
 \etmz
 \vspace{-2ex}
 \qed

\bsk

In a similar way, we also obtain the following non-approximability
 result concerning $\tau_2$.

\bc
 The value\/ $\tau_2(\cH)$
  does not have a polynomial-time\/
  $((1-\eps)\ln m)$-approximation on hyperstars,
   unless {\sf  NP}\,$\subseteq$\,{\sf DTIME}$(n^{\cO(\log \log
   n)})$.
 \ec

 \pf
By Proposition \ref{star-decrem}, the approximation of $\tau_2(\cH)$
 on hyperstars $\cH$ is as hard as that of $\tau(\cH^-)$ on
 general hypergraphs $\cH^-$.
 \qed

\brm
 In connection with Theorem \ref{ratio} one may observe that,
 even if we restrict the problem instances to 3-uniform hypergraphs
 in which each vertex pair is contained in at most three edges,\/
 $\UU(\cH)$ does not admit a PTAS.
 This follows from the fact that the determination of\/ $\aaa(G)$
 is\/ \mbox{{\sf MAX SNP}}-complete on graphs of maximum degree 3,
 by the theorem of Berman and Fujito \cite{BF}.
\erm

\section{Concluding remarks}

Our results on hyperstars show that $\dec(\cH)$ admits a much better
approximation than $\UU(\cH)$ does. In a way this fact is in analogy
with the following similar phenomenon in graph theory: The
independence number\/ $\aaa(G)$ is not approximable within\/
$n^{1-\eps}$, but\/ $\tau(G)=n-\aaa(G)$ admits a polynomial-time
2-approximation because\/ $\nu(G)\le\tau(G)\le 2\nu(G)$, and the
matching number\/ $\nu(G)$ can be determined in polynomial time.
In this way, both comparisons $\dec(\cH)$ with $\UU(\cH)$
and $\tau(G)$ with $\aaa(G)$ demonstrate that there
 can occur substantial difference between the approximability of
a graph invariant and its complement.

Perhaps hypertrees with not very large edges admit some
 fairly efficient algorithms:

\bop
 Determine the largest integer\/ $r$ such that there is a PTAS
  to approximate the value of\/ $\UU(\cH)$
   for hypergraphs\/ $\cH$ in which every edge has
  at most\/ $r$ vertices.
\eop

Our results imply that $r\le 6$ is necessary. From below, a very
easy observation shows that for $r=2$ there is a linear-time
 algorithm, because for graphs $G$, the value
of $\UU(G)$ is precisely the number of connected components.

For hypertrees with non-restricted edge size, the following open
question seems to be the most important one:

\bop
 Is there a polynomial-time\/ $o(n)$-approximation for\/
  $\UU$ on hypertrees?
\eop

\end{document}